\RequirePackage{amsmath}
\RequirePackage{fix-cm}
\documentclass{svjour3}                     
\pdfoutput=1
\usepackage{authblk}
\usepackage{lineno,hyperref}
\usepackage{graphicx}
\usepackage{epstopdf}
\usepackage{bm}
\usepackage[round]{natbib}
\usepackage{lineno}

\DeclareGraphicsExtensions{.png,.jpeg,.pdf}
\numberwithin{equation}{section}
\allowdisplaybreaks
 
\newcommand*\patchAmsMathEnvironmentForLineno[1]{%
  \expandafter\let\csname old#1\expandafter\endcsname\csname #1\endcsname
  \expandafter\let\csname oldend#1\expandafter\endcsname\csname end#1\endcsname
  \renewenvironment{#1}%
     {\linenomath\csname old#1\endcsname}%
     {\csname oldend#1\endcsname\endlinenomath}}%
\newcommand*\patchBothAmsMathEnvironmentsForLineno[1]{%
  \patchAmsMathEnvironmentForLineno{#1}%
  \patchAmsMathEnvironmentForLineno{#1*}}%
\AtBeginDocument{%
\patchBothAmsMathEnvironmentsForLineno{equation}%
\patchBothAmsMathEnvironmentsForLineno{align}%
\patchBothAmsMathEnvironmentsForLineno{flalign}%
\patchBothAmsMathEnvironmentsForLineno{alignat}%
\patchBothAmsMathEnvironmentsForLineno{gather}%
\patchBothAmsMathEnvironmentsForLineno{multline}%
}
 
\journalname{Noname}
 
\begin{document}
 
\title{The Role of Stochasticity in Noise-Induced Tipping Point Cascades: A Master Equation Approach}
 
\titlerunning{The Role of Stochasticity in Noise-Induced Tipping Point Cascades}
 
\author{Abhishek Mallela \and
        Alan Hastings
}
 
\institute{Abhishek Mallela \at
              Department of Mathematics, University of California Davis, Davis, CA 95616 \newline
              \email{amallela@ucdavis.edu}           
                   \and
                    Alan Hastings \at
                    Department of Environmental Science and Policy, University of California Davis, Davis, CA 95616
                    \at Santa Fe Institute, Santa Fe, NM 87501
}
 
\date{}
 
\maketitle
 
\abstract
Tipping points have been shown to be ubiquitous, both in models and empirically in a range of physical and biological systems. The question of how tipping points cascade through systems has been less well studied and is an important one. A study of noise-induced tipping, in particular, could provide key insights into tipping cascades. Here, we consider a specific example of a simple model system that could have cascading tipping points. This model consists of two interacting populations with underlying Allee effects and stochastic dynamics, in separate patches connected by dispersal, which can generate bistability. From an ecological standpoint, we look for rescue effects whereby one population can prevent the collapse of a second population. As a way to investigate the stochastic dynamics, we use an individual-based modeling approach rooted in chemical reaction network theory. Then, using continuous-time Markov chains and the theory of first passage times, we essentially approximate, or emulate, the original high-dimensional model by a Markov chain with just four states, where each state corresponds to a combination of population thresholds. Analysis of this reduced model shows when the system is likely to recover, as well as when tipping cascades through the whole system.

\section {Introduction}
Many systems in nature can transition into a qualitatively distinct dynamic regime when a critical threshold is approached. The associated threshold of such a system, defined in the context of bifurcation theory, is the bifurcation point or tipping point of the system. Tipping points manifest in systems including lake eutrophication (ecology), social contagion (sociology), disease spread (epidemiology), epileptic seizures (physiology), stock market crashes (finance), and even the earth system (climate science) \citep{scheffer, scheffer-2012, oregan-2018, lenton-2020, klose-2020}. Tipping points  arise in the presence of strongly self-amplifying (mathematically positive) feedbacks \citep{lenton-2020}. This characterization of tipping points is important in the sense that only sufficiently strong, self-propelling feedbacks are recognized.

A prime example of a tipping point is an ecological system with an Allee effect \citep{courchamp, etc, johnson-2018, vortkamp-2020}. Allee effects occur in populations with low abundances and are believed to be commonplace in ecological systems \citep{courchamp, drake-2011}. In a population exhibiting an Allee effect, the per capita growth rate is a unimodal function of the population abundance with a global maximum. The sign of the growth rate at low population levels distinguishes a weak Allee effect from a strong one. In particular, a weak Allee effect does not result in a negative growth rate for small population sizes but its strong counterpart presents with a negative rate. Therefore, the existence of a strong Allee effect implies the existence of a critical threshold for survival. The possibility of alternative stable states in systems that are analogous to those with an Allee effect has significant implications at many levels, from the microscopic scale of budding yeast \citep{dai-2012} to the macroscopic scale of tipping elements for the Earth system \citep{lenton-2020, klose-2020}.

In order to investigate tipping points induced by sources of variation in systems with an Allee effect, one needs to formulate a model that can exhibit either a weak or strong Allee effect. As the quality of the environment is degraded, a weak Allee effect becomes a strong one through a transcritical bifurcation. A fundamental question in this setting concerns the propagation of tipping points through an ecosystem consisting of multiple patches. To what extent are interacting populations interdependent and how is this relationship influenced by the parameters governing the model behavior? Stated more simply, how do tipping points cascade through systems?

A simple analysis considers multiple populations in a network that are connected by passive, symmetric diffusion. In this work, we study a model consisting of two populations to analyze how their dynamics are related, conditional on the quality of their internal environments. One case of interest is when both populations are above the Allee threshold to see when the collapse of either population is less likely than for isolated systems.  A second case of interest is when one population has fallen below the Allee threshold and we ask whether the next transition puts both populations above, or alternatively, below, the Allee threshold. A similar study was conducted in a deterministic setting \citep{johnson-2018}. That study laid the foundations for an eventual treatment in a stochastic setting. Hence, in our work, we adopt the definition of ecological resilience formulated by Holling \citep{holling-1973} as well as the framework proposed by Johnson and Hastings, but we also account for stochasticity in the system.

The literature on single, isolated tipping points is vast. However, studies of cascades of tipping points are less common. We begin with the background needed to formulate a stochastic, individual-based model that accounts for Allee effects of varying strengths. Our model follows a Markovian birth-death process, rooted in the theory of chemical reaction networks. We find that this representation naturally lends itself to a treatment with Markov chains and first passage times. Then, we employ dimensionality reduction techniques and various approximations to derive a reduced model that we study in detail. The section that follows consists of the results obtained using the above methods. Finally, we summarize our findings with conclusions.

\section {Model Description}
Much of the deterministic modeling work done on systems with one population that manifest Allee effects are based on either phenomenology or empirical observations. The crux of these models dates back to 1954 \citep{odum-1954}, where the observed per capita growth rate was fit using a suitable function. The general form of a deterministic model in continuous time is the following ordinary differential equation:
\begin{equation}
\dfrac{d \rho(t)}{dt} = \rho f(\rho),
\end{equation}
where $\rho(t)$ denotes the average number of individuals at time $t,$ and $f(\rho)$ is a function specifying the form of the per-capita population growth rate at size $\rho$. There have been many functional forms proposed for $f(\rho)$ in the literature. A review of various specifications used in models can be found \citep{boukal-2002}.

In particular, a simple and important model was proposed early in the twentieth century \citep{volterra-1938}, where $f(\rho)$ is a quadratic polynomial function of $\rho$. The underlying assumptions of the model are as follows. Given a constant sex ratio, the number of meetings between males and females is proportional to $\rho^2$. The ratio of births to meetings can be affected by the population density and is hence assumed to be linearly decreasing in $\rho$. Also, the influx and efflux of individuals are represented by birth and death events that occur at constant per capita rates. Hence, the Volterra model takes the following form:

\begin{equation}
\dfrac{d \rho(t)}{dt} = -a_1 \rho + (a_2-a_3\rho)\rho^2 = -a_1 \rho + a_2 \rho^2 - a_3 \rho^3,
\end{equation}
where $a_2, a_3 > 0$. The sign of $a_1,$ however, is determined by the magnitude of the difference between the birth and death rates. If we define the two real-valued roots,
\begin{subequations}
\begin{align}
k_1 &= \dfrac{1}{2a_3}[a_2 - \sqrt{a_2^2 - 4a_1a_3}], \\
k_2 &= \dfrac{1}{2a_3}[a_2 + \sqrt{a_2^2 - 4a_1a_3}],
\end{align}
\end{subequations}
with $a_2^2 > 4a_1a_3$, then the model is often shown in the following form \citep{courchamp-1999}
\begin{equation}
\label{base}
\dfrac{d \rho(t)}{dt} = a_1\rho\left(1-\dfrac{\rho}{k_2}\right)\left(\dfrac{\rho}{k_1} - 1\right),
\end{equation}
resembling the logistic model with the addition of a new unstable steady state, $\rho = k_1.$

If $a_1 > 0$, the Volterra model has three steady states:
two stable states at $\rho = 0$ and $\rho = k_2$ and an unstable state at $\rho = k_1$. In this case, if the initial population size exceeds $k_1,$ the population grows over time and converges to the stable steady state $\rho = k_2$, the carrying capacity of the system. If the initial population size is less than $k_1$, the population decays over time to the stable steady state $\rho = 0$ and goes extinct. This scenario describes the strong Allee effect. The dynamics when $a_1 < 0$ are different. In this scenario, $k_1$ becomes negative and the steady state $\rho = k_1$ lacks biological meaning. Here, the Volterra model has only two steady states, an unstable state at $\rho = 0$ and a stable state at $\rho = k_2$. This case corresponds to the weak Allee effect.

In order to investigate the effects of internal fluctuations or demographic stochasticity in a two-population system with both weak and strong Allee effects, we consider the temporal evolution of the system as specified by a Markovian birth-death process. Demographic stochasticity is included both in the dynamics of the locations and in the dispersal parameter. We account for the individual reaction kinetics explicitly, in a mechanistic manner, without relying on phenomenological considerations. Here, we follow the approach of M{\'e}ndez and colleagues \citep{mendez-2019} by casting the system as a chemical reaction network that results in an individual-based model (IBM). The minimal IBM that displays both the weak and strong Allee effect and also accounts for dispersal can be described as follows. It consists of two birth processes (linear and binary birth), a ternary competition process, a linear death process, and an exchange process. We provide our reaction scheme below.

\begin{subequations}
\label{reaction-scheme}
\begin{align} 
&X_1 \overset{\mu_1}{\rightharpoonup} (1+b)X_1 \\
2&X_1 \overset{\lambda_1}{\rightharpoonup} (2+a)X_1 \\
&X_1 \overset{\gamma_1}{\rightharpoonup} \emptyset \\
3&X_1 \overset{\tau_1}{\rightharpoonup} (3-c)X_1 \\
&X_1 \overset{d}{\rightharpoonup} X_2 \\
&X_2 \overset{\mu_2}{\rightharpoonup} (1+b)X_2 \\
2&X_2 \overset{\lambda_2}{\rightharpoonup} (2+a)X_2 \\
&X_2 \overset{\gamma_2}{\rightharpoonup} \emptyset \\
3&X_2 \overset{\tau_2}{\rightharpoonup} (3-c)X_2 \\
&X_2 \overset{d}{\rightharpoonup} X_1
\end{align}
\end{subequations}

The first reaction is a linear birth process, which occurs at a constant rate $\mu_1$, and describes the baseline reproductive success of the first population in the absence of cooperative effects. It accounts for the fact that the typical individual produces $b$ offspring that reach reproductive age. The second reaction is a binary process that occurs at a constant rate $\lambda_1$. It describes cooperative interactions, such as breeding, antipredator behavior, or foraging, that result in producing $a$ additional offspring which reach reproductive age. The third reaction is a linear death process, occurring at constant rate $\gamma_1$, which accounts for mortality due to natural causes. The fourth reaction is a ternary competition process, accounting for the results of overcrowding and resource depletion, where $c$ individuals die at rate $\tau_1$. Note that $c = 1, 2, 3$ are the only meaningful values. The next reaction is an exchange process of symmetric dispersal between the two populations. This occurs at a constant rate of $d$. The last five reactions in the scheme describe the dynamics of the second population, respectively.

The reaction scheme \eqref{reaction-scheme} defines a Markovian process, and the temporal evolution of $P(n_1,n_2,t)$, the probability of having $n_i$ individuals from the $i^{th}$ population at time $t$ for $i = 1,2$, is described by the following master equation, also known as the forward Kolmogorov equation \citep{gardiner}:

\begin{equation}
\label{master}
\dfrac{dP(\bm{n}, t)}{dt} = \sum\limits_{\bm{r}} [W(\bm{n} - \bm{r}, \bm{r}) P(\bm{n} - \bm{r}, t) - W(\bm{n}, \bm{r}) P(\bm{n}, t)],
\end{equation}
where $P(\bm{n} < \bm{0}, t) = P(n_1 < 0, n_2 < 0, t) = 0.$ Here $W(\bm{n},\bm{r})$ are the transition rates between the states with $\bm{n}$ and $\bm{n} + \bm{r}$ individuals, where $$\bm{r} = \{\bm{r_1},\bm{r_2},\ldots,\bm{r_{10}}\} = \{(b,0),(a,0),(-1,0),(c,0),(-1,0),(0,b),(0,a),(0,-1),(0,c),(0,-1)\}$$ is the vector of transition increments corresponding to the system given by Eq. \eqref{reaction-scheme}. The transition rates corresponding to each reaction, $W(\bm{n},\bm{r})$, are obtained from the reaction kinetics (\cite{vankampen,gardiner}):

\begin{subequations}
\label{transition-rates}
\begin{align}
W(\bm{n},\bm{r_1}) &= \mu_1 n_1, \\
W(\bm{n},\bm{r_2}) &= \dfrac{\lambda_1}{2} n_1 (n_1-1), \\
W(\bm{n},\bm{r_3}) &= \gamma_1 n_1, \\
W(\bm{n},\bm{r_4}) &= \dfrac{\tau_1}{6} n_1 (n_1-1) (n_1-2), \\
W(\bm{n},\bm{r_5}) &= d n_2, \\
W(\bm{n},\bm{r_6}) &= \mu_2 n_2, \\
W(\bm{n},\bm{r_7}) &= \dfrac{\lambda_2}{2} n_2 (n_2-1), \\
W(\bm{n},\bm{r_8}) &= \gamma_2 n_2, \\
W(\bm{n},\bm{r_9}) &= \dfrac{\tau_2}{6} n_2 (n_2-1) (n_2-2), \\
W(\bm{n},\bm{r_{10}}) &= d n_1.
\end{align}
\end{subequations}

Deterministic ODEs for the average population size can be obtained from Eq. \eqref{master}. Multiplying Eq. \eqref{master} by $n_1 n_2$, using transition rates \eqref{transition-rates}, and summing over all values of $n_1$ and $n_2$, we find 
\begin{equation}
\label{mean-field}
\dfrac{d \rho_i}{dt} = (\mu_i b - \gamma_i - d)\rho_i + \dfrac{a \gamma_i}{2}\rho_i^2 - \dfrac{c \tau_i}{6} \rho_i^3,
\end{equation}
for $i = 1, 2$, where $\rho_i = \langle n_i \rangle$ is the mean number of individuals in population $i$. We note that this deterministic equation holds strictly when the demographic fluctuations vanish, which occurs in the thermodynamic limit as the population sizes increase to infinity. Writing $\bm{k_1} =(k_{11}, k_{12})$ and $\bm{k_2} =(k_{21}, k_{22})$, Eq. \eqref{mean-field} can be cast in the form of Eq. \eqref{base} with the definitions
\begin{align*}
\bm{k_1} &= \dfrac{3}{2c\tau_i}[a \lambda_i - \sqrt{a^2 \lambda_i^2 +8c\tau_i(b\mu_i - \gamma_i-d)/3}], \\
\bm{k_2} &= \dfrac{3}{2c\tau_i}[a \lambda_i + \sqrt{a^2 \lambda_i^2 +8c\tau_i(b\mu_i - \gamma_i-d)/3}]
\end{align*}
At the deterministic level, the enzymatic reaction scheme \eqref{reaction-scheme} of the IBM gives rise to the Volterra rate equation Eq. \eqref{base}.

For the sake of simplicity, in what follows, we focus on the simplest version of this IBM. Namely, we treat the Markovian process as a single-step process with $a = b = c = 1.$ The set of reactions \eqref{reaction-scheme} then becomes

\begin{subequations}
\label{reaction-scheme-new}
\begin{align} 
&X_1 \overset{\mu_1}{\rightharpoonup} 2X_1 \\
2&X_1 \overset{\lambda_1}{\rightharpoonup} 3X_1 \\
&X_1 \overset{\gamma_1}{\rightharpoonup} \emptyset \\
3&X_1 \overset{\tau_1}{\rightharpoonup} 2X_1 \\
&X_1 \overset{d}{\rightharpoonup} X_2 \\
&X_2 \overset{\mu_2}{\rightharpoonup} 2X_2 \\
2&X_2 \overset{\lambda_2}{\rightharpoonup} 3X_2 \\
&X_2 \overset{\gamma_2}{\rightharpoonup} \emptyset \\
3&X_2 \overset{\tau_2}{\rightharpoonup} 2X_2 \\
&X_2 \overset{d}{\rightharpoonup} X_1
\end{align}
\end{subequations}

The mean-field rate equation corresponding to \eqref{reaction-scheme-new} is 
\begin{equation}
\label{mean-field-new}
\dfrac{d \rho_i}{dt} = (\mu_i - \gamma_i - d)\rho_i + \dfrac{\gamma_i}{2}\rho_i^2 - \dfrac{\tau_i}{6} \rho_i^3,
\end{equation}
for $i=1, 2$. \\

For this set of reactions, the master equation can be explicitly obtained as
\begin{align}
\label{master-new}
\dfrac{dP(n_1,n_2,t)}{dt} 
&= (n_1-1)\left[\dfrac{\lambda_1}{2}(n_1-2) + \mu_1\right] P(n_1-1,n_2,t) \\ \notag
&+ (n_1+1)\left[\dfrac{\tau_1}{6}n_1(n_1-1) + \gamma_1+d\right] P(n_1+1,n_2,t) \\ \notag
&+ (n_2-1)\left[\dfrac{\lambda_2}{2}(n_2-2) + \mu_2\right] P(n_1,n_2-1,t) \\ \notag 
&+ (n_2+1)\left[\dfrac{\tau_2}{6}n_2(n_2-1) + \gamma_2+d\right] P(n_1,n_2+1,t) \\ \notag 
&- n_1\left[\dfrac{\lambda_1}{2}(n_1-1)+\mu_1+\dfrac{\tau_1}{6}(n_1-1)(n_1-2)+\gamma_1+d\right]P(n_1,n_2,t) \\ \notag 
&-n_2\left[\dfrac{\lambda_2}{2}(n_2-1)+\mu_2+\dfrac{\tau_2}{6}(n_2-1)(n_2-2)+\gamma_2+d\right]P(n_1,n_2,t)
\end{align}

We note that the master equation \eqref{master-new} includes only single-step processes where the transitions take place between the states $\bm{n}$ and $\bm{n} \pm (1,0)$ or $\bm{n}$ and $\bm{n} \pm (0,1)$. We can define new dimensionless quantities in terms of the reaction rates as follows:
\begin{equation}
\label{dimensionless}
N_i = \dfrac{3 \lambda_i}{2 \tau_i}, \ \ \ \delta_i^2 = 1 + \dfrac{8\tau_i(\mu_i-\gamma_i-d)}{3\lambda_i^2}, \ \ \ \ R_0^{(i)} = \dfrac{\mu_i}{\gamma_i+d},
\end{equation}

Note that $N_i$ defines the scale of the typical size of population $i$ prior to extinction. The identities \eqref{dimensionless} establish a relation between the microscopic $(\lambda_i, \mu_i, \gamma_i, \tau_i, d)$ and macroscopic $(N_i, \delta_i, R_0^{(i)})$ parameters, which are obtainable through field observations. We now realize that the IBM displays both types of Allee effects:
\begin{subequations}
\begin{align}
&\text{Weak Allee}: \ \ \ \mu_i > \gamma_i + d \ \ \text{or} \ \ R_0^{(i)} > 1, \\
&\text{Strong Allee}: \ \ \ \mu_i < \gamma_i + d \ \ \text{or} \ \ R_0^{(i)} < 1.
\end{align}
\end{subequations}
Note that for $R_0^{(i)} > 1 \ (R_0^{(i)} < 1)$ we have $\delta_i > 1 \ (\delta_i <1)$, where
for the strong Allee effect, we must also demand that $\delta_i > 0.$ \\

As mentioned in the Introduction, our starting point for the model used in this work is the deterministic skeleton for the non-dimensionalized model described by Johnson and Hastings, reproduced here for convenience.
\begin{subequations}
\begin{align}
\label{johnson1}
\dot{X_1} = X_1(\beta_1-(X_1-1)^2 )+D(X_2-X_1) \\
\label{johnson2}
\dot{X_2} = X_2 (\beta_2-(X_2-1)^2 )+D(X_1-X_2)
\end{align}
\end{subequations}

In the model above, the parameter $\beta_i$ represents a measure for the quality of the environment by the population denoted by $X_i$. The parameter $D$ denotes passive diffusion in the system and is an indicator of network connectivity. See \citep{johnson-2018} for a detailed exposition of the model.

In order to make our subsequent analyses feasible, we aim to reduce the dimensionality of the parameter space for our stochastic model. So, for the sake of illustration, we treat the stochastic rate parameters in \eqref{reaction-scheme-new} as identical for both populations. This assumption appears to be reasonable since a compelling justification for heterogeneity in the population parameters is unlikely. Moreover, this assumption is consistent with the model parameterization in the work by Johnson and Hastings.

Thus, our matching scheme can be written as follows:
\begin{subequations}
\begin{align}
\label{matching-scheme}
\mu_1 &:= \beta_1, \\
\mu_2 &:= \beta_2, \\
\lambda &:= \lambda_1 = \lambda_2 \equiv 4, \\
\gamma &:= \gamma_1 = \gamma_2 \equiv 1, \\
\tau &:= \tau_1 = \tau_2 \equiv 6, \\
d &:= D, \\
\tilde{N}  &:= \tilde{N}_1 = \tilde{N}_2
\end{align}
\end{subequations}

So for $i = 1$ and $2$,
\begin{equation}
\label{dimensionless-new}
\tilde{N} = 1, \ \ \ \delta_i^2 = \beta_i-D, \ \ \ \ R_0^{(i)} = \dfrac{\beta_i}{D + 1}.
\end{equation}

In what follows, we are interested in the case of bistability, which is manifest in the case of the strong Allee effect. Note that this necessitates the following condition:
\begin{equation}
D < \beta_i < D+1
\end{equation}

\section{Methods}
We begin by noting that our stochastic model operates over a two-dimensional state space. We argue that the dimension of the state space is as low as possible but nevertheless captures the desired phenomena. Due to the presence of Allee effects, there are cubic nonlinearities in the system. So, the existence of an analytic solution to the system is highly unlikely. To circumvent this issue, Johnson and Hastings conducted numerical simulations using the deterministic skeleton in order to produce bifurcation diagrams. They focused on how the transition rate between the two patches, or $D$, determine the bifurcation structure for this system. In this study, however, we restrict our attention to the case of bistability, which is properly addressed with a stochastic model.

As discussed in the previous section, we chose to adapt a master equation approach for our model \citep{mendez-2019}. Instead of specifying a model with carrying capacities, we used an individual-based modeling approach using a chemical reaction network. This allows for a fine-grained representation of the underlying discrete, stochastic process. Using this approach, we wrote down the two-dimensional chemical master equation for the process.

Given that this stochastic process is a continuous-time Markov chain (CTMC), it can be explicitly described by a generator $\mathcal{Q}$-matrix with a countable state space. A nice feature of most ecological models is that they are built around processes that will approach a compact set exponentially fast. Any reasonable ecological model should not have unbounded population growth. Density dependence in ecology models is typically what gives this. This implies that the model effectively operates over a finite state space as the probability of arbitrarily large populations is negligibly small. Using this insight, we were able to obtain a finite-state CTMC in two dimensions.

Since the multi-dimensional master equation was relatively unwiedly to work with, we sought to reduce the two-dimensional state space to one dimension \citep{allen-2010, allen-2003}. Denoting $N := N_1 = N_2$ as the maximum number of individuals in either population, the specific mapping function used was
\begin{equation}
\label{mapping}
f\left((x,y)\right) = (N+1)x + y + 1, \ \ \ \ \ 0 \le x,y \le N.
\end{equation}
We could also exploit the sparsity of the banded $\mathcal{Q}$-matrix to drastically simplify the effort needed in computations \citep{doorn-2013}. Specifically, the $\mathcal{Q}$-matrix has size $(N+1)^2 \times (N+1)^2$ with $7 N^2 + 4N - 2$ nonzero entries, yielding a matrix density of $\mathcal{O}(N^{-2}).$ Thus, the sparsity of the matrix increases quadratically with $N.$

For the ensuing analyses, we obtained population thresholds as follows. The high thresholds were varied from $H := H_1 = H_2 = 2$ to $H = N$. Here, $N$ is defined as before, and can be understood as a system size parameter. The low thresholds were varied from $L := L_1 = L_2 = \tilde{N} = 1$ to $L = H-1$. This ensures that $L$ is always less than $H.$ Note that the smallest $L$ threshold, $L = 1,$ is representative of quasi-extinction, or the typical size of either population before extinction.

In order to probe the system under consideration, we computed the mean first passage times (MFPTs) for the model with the state-space parameterized by $N$ for all combinations of the parameters $\beta_i, D, L,$ and $H$ \citep{chou-2014, polizzi-2016}. This was done as follows. First, we formed the $\mathcal{Q}$-matrix for each point in parameter space. Then, the rows and columns corresponding to the trap states (e.g. extinction at $(x = 0, y = 0)$) were removed. We also formed a vector of initial state probabilities $p_0$ governing the subsequent evolution of the CTMC. This vector has $(N+1)^2$ entries. The corresponding entry was removed from this vector. Next, using the resulting truncated $\mathcal{Q}$-matrix, $\tilde{\mathcal{Q}}$, we computed the matrix of mean residence times, or  $-\tilde{\mathcal{Q}}^{-1}$. Finally, we computed the sum of the entries in $-\tilde{\mathcal{Q}}^{-1} p_0$ to yield the MFPT from the initial state to the desired trap state.

We could then construct a compartmental system with a reduced state space consisting of just four states: $HH$, $HL$, $LL$, and $LH$. Each of these states corresponds to a combination of population thresholds. For instance, $HL$ means that the first population is at a high abundance and the second population is at a low level. The MFPTs from the original model were used as input rates for the transition rate $\mathcal{S}$-matrix of the reduced model. Thus we used an emulator, or meta-model, as a proxy to analyze the original system.

\section{Results}
\begin{figure}
\centering
\includegraphics[width=\linewidth]{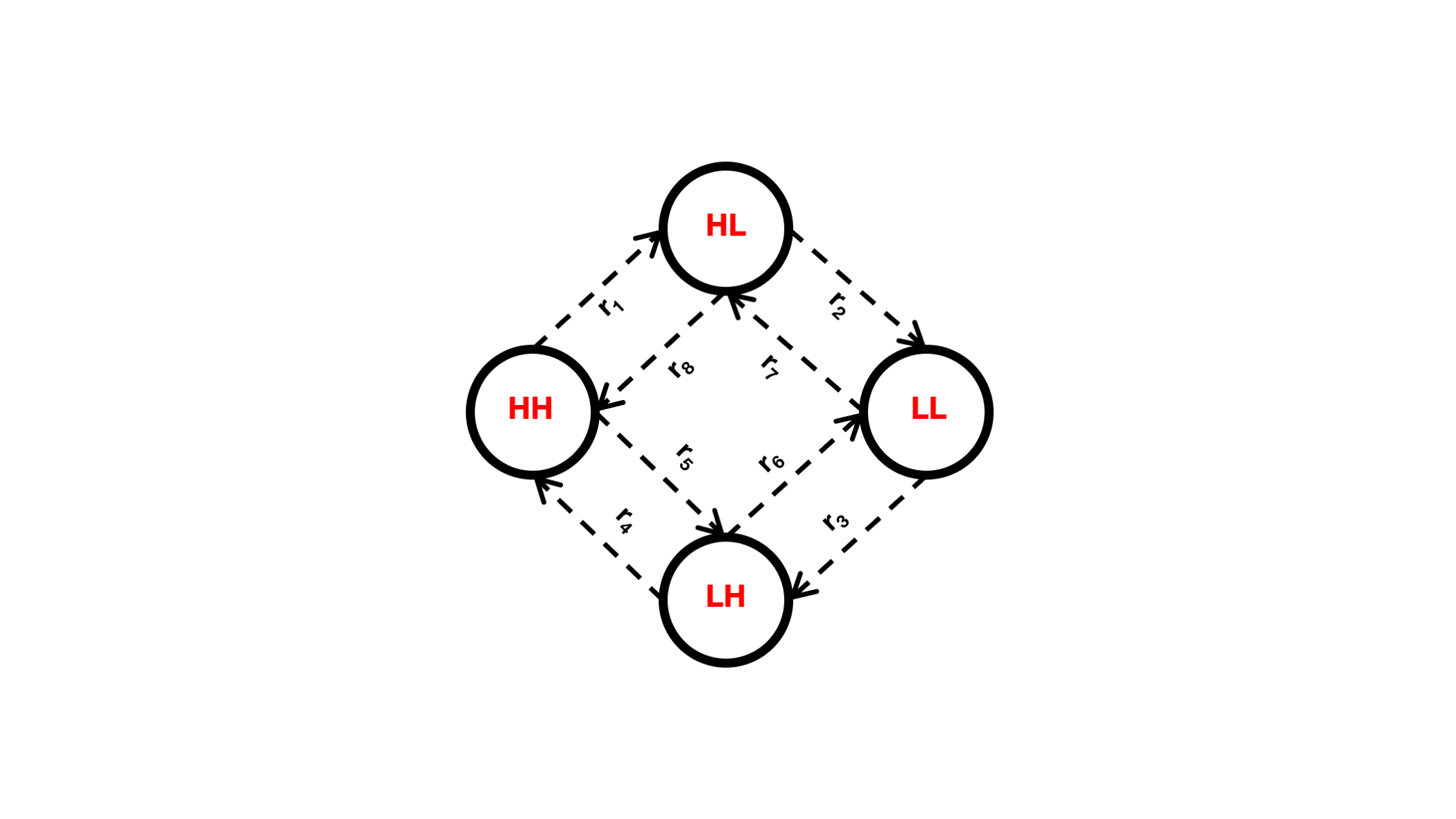}
\caption{\textbf{Schematic diagram of the complete emulator.} Each state of the compartmental system is shown as a black circle and denoted by the letters in red. The letter $H$ corresponds to the high Allee threshold and $L$ corresponds to the low threshold. The first (second) letter refers to the first (second) population. Each $r_i$ for $i = 1, \ldots, 8$ describes the transition rate, or inverse mean first passage time, between appropriate states.}
\label{one}
\end{figure}

Throughout this work, we analyze the simplest possible system capable of exhibiting noise-induced tipping point cascades in the vicinity of saddle-node bifurcations. We can use the emulator described previously to construct the schematic diagram in Fig. \ref{one}. In the diagram, each $r_i$ for $i = 1, \ldots, 8$ represents the rate of the transition between the relevant compartments. Each rate can be computed as the inverse of its corresponding mean first passage time. Given that the process is represented as a CTMC, we can write down the transition rate matrix $\mathcal{S}$ for the emulator. The $n$-step transition probability matrix $(\mathcal{P}^n)_{ij}$ gives the $n$-step transition probabilities of the system, from the state corresponding to the $i^{th}$ row to the state corresponding to the $j^{th}$ column. \\

The $\mathcal{S}$-matrix is given as:
\begin{equation*}
\mathcal{S} = 
\begin{pmatrix}
-r_1-r_5 & r_1 & r_5 & 0 \\
r_8 & -r_2-r_8 & 0 & r_2 \\
r_4 & 0 & -r_4-r_6 & r_6 \\
0 & r_7 & r_3 & -r_3-r_7
\end{pmatrix}
\end{equation*}
where the ordering of the states is $(HH, HL, LH, LL)$. Note that $\mathcal{S}$ satisfies the properties of a generator matrix, namely:
\begin{itemize}
\item All off-diagonal elements are non-negative.
\item All diagonal elements are negative.
\item Each row sums to zero.
\end{itemize}

\begin{figure}
\centering
\includegraphics[width=\linewidth]{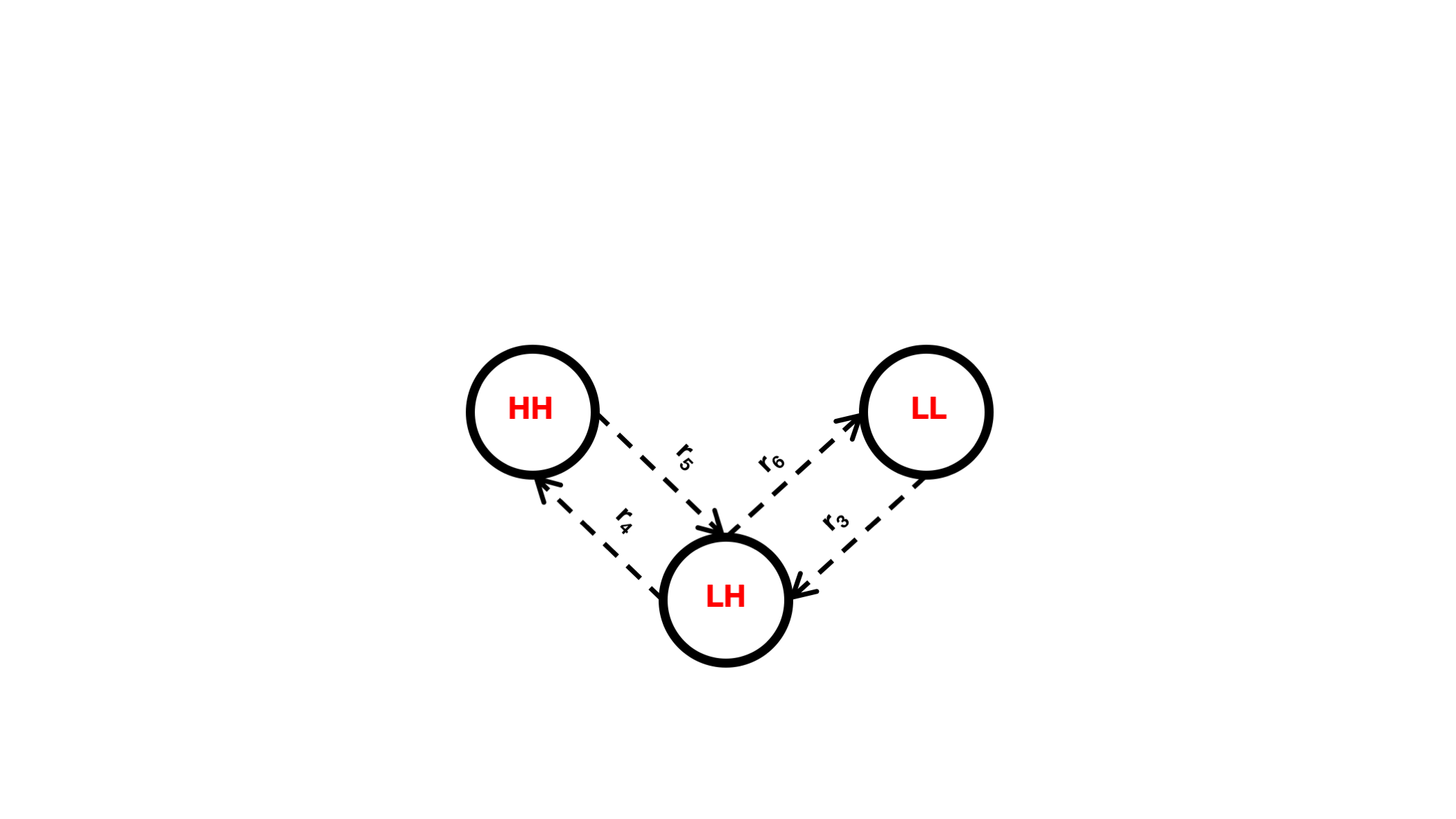}
\caption{\textbf{Schematic diagram of the reduced emulator.} Identical diagram as in Fig. \ref{one}, except for the omission of the $HL$ state. See the text for justification.}
\label{two}
\end{figure}

Since the two populations $X_1$ and $X_2$ interact with each other via the symmetric dispersal parameter $D$ and the thresholds for both populations are equivalent, we can treat the states $LH$ and $HL$ as equivalent. Hence, we can simplify our formalism by omitting the $HL$ state (Fig. \ref{two}):

\begin{equation*}
\mathcal{\tilde{S}} = 
\begin{pmatrix}
-r_5 & r_5 & 0 \\
r_4 & -r_4-r_6 & r_6 \\
0 & r_3 & -r_3
\end{pmatrix}
\end{equation*}
where the ordering of the states is now $(HH, LH, LL)$. \\

Now, we obtain the following results:
\begin{align}
&P(HH \to LH) = 1 \ \text{with} \ MFPT = \frac{1}{r_5} \label{41}\\
&P(LH \to HH) = \frac{r_4}{r_4+r_6} \ \text{with} \ MFPT = \frac{1}{r_4} \label{42}\\
&P(LH \to LL) = \frac{r_6}{r_4+r_6} \ \text{with} \ MFPT = \frac{1}{r_6} \label{43}
\end{align}

By the multiplicative rule for probabilities, we obtain:
\begin{align}
&P(HH \to LL) = \frac{r_6}{r_4+r_6} \ \text{with} \ MFPT = \frac{1}{r_5} + \frac{1}{r_6} \\
&P(HH \to HH) = \frac{r_4}{r_4+r_6} \ \text{with} \ MFPT = \frac{1}{r_5} + \frac{1}{r_4}
\end{align}

Hence, the probability of a system collapse differs from the probability of system resiliency in a fundamental manner that is quantified by the parameters 
\begin{align}
r_6 = \frac{1}{MFPT(LH \to LL)} , \ \ r_4 = \frac{1}{MFPT(LH \to HH)}
\end{align}

If $r_4 = r_6$, then $P(HH \to HH) = r_4/(r_4+r_6) = r_6/(r_4+r_6) = P(HH \to LL) =1/2$ so the system is equally likely to collapse or recover. We observe the following chain of implications in the case that $r_4 > r_6$. The likelihood of recovery is larger than the threshold value of $\frac{1}{2}$, which implies that the MFPT decreases in the vicinity of the tipping point of the system. This in turn implies that $r_4 > r_6$. \\

Our aim henceforth is to explore the multi-dimensional region of parameter space that corresponds to the (desired) probability of system resiliency. In symbols, for a fixed $N$, we want to identify all combinations $(D,\beta_1,\beta_2,H, L)$ where $D >0, \ \beta_1 \in (D, D+1), \ \beta_2 \in (D, D+1),$ $H \in [2,N],$ and $L \in [1,H-1]$ and the system recovers with the aid of rescue effects. In order to perform these numerical simulations, we began by specifying an upper bound of $1$ for the dispersal parameter $D.$ This value for $D$ could be considered as large, since $D$ should be commensurate with the per-capita rate of the system, which was chosen as $r = 1$ in the non-dimensionalization of the model Eqs. \eqref{johnson1}-\eqref{johnson2} \citep{johnson-2018}. This implies that $D \in (0, 1).$ Then, we chose $N = 10.$ This was the largest $N$ such that the largest condition number among the space of all matrices $\tilde{\mathcal{Q}}$, taken over the set of simulated parameter values, was less than $10^{7}$. For our chosen value of $N,$ we found this condition number to be $8.56 \times 10^{6}.$ We note that numerical linear algebra was used here instead of Monte-Carlo methods, including the Doob-Gillespie algorithm \citep{gillespie-1976, gillespie-1977}). The primary justification for this choice of method was that some of the simulated mean first passage times were very large, causing prohibitively long runtimes with Monte-Carlo simulations. Thus, $\beta_1$ and $\beta_2$ were restricted accordingly between $D$ and $D+1$. Finally, the high and low thresholds were varied with unit spacing. All computations were performed in MATLAB \citep{MATLAB:2020}. \\

Defining $r = P(LH \to HH) = \frac{r_4}{r_4+r_6},$ we can analyze the parameter space. For the sake of illustration, we chose three values for the habitability parameters, for each value of $D$ (i.e. $\beta_i = D+0.01, D+0.5, D+0.99$). Upon analysis of the simulated data set, we found that the minimum value of $r$ over the entire parameter space was approximately $0.275$, and the maximum value was $1$. This indicates that it is possible to guarantee system recovery for some set of parameter values, but there is no set of parameter values that ensures system failure. In other words, the system is never more than $72.5\%$ likely to fail. Next, we varied the threshold value $\eta$, such that $r > \eta$, for $\eta = 0.9, 0.95, 0.99$. How weak can the environment be, or equivalently, which resource constraints will not obstruct our goal of achieving $r > \eta$? We can intuit that the percentage of parameter space for which the system recovers will decrease monotonically with the threshold $\eta$. Naturally, we would like to find a region of parameter space such that $\eta$ is large and the proportion of parameter space $\nu$ for which $r > \eta$ is large. \\

\begin{figure}
\centering
\includegraphics[width=\linewidth]{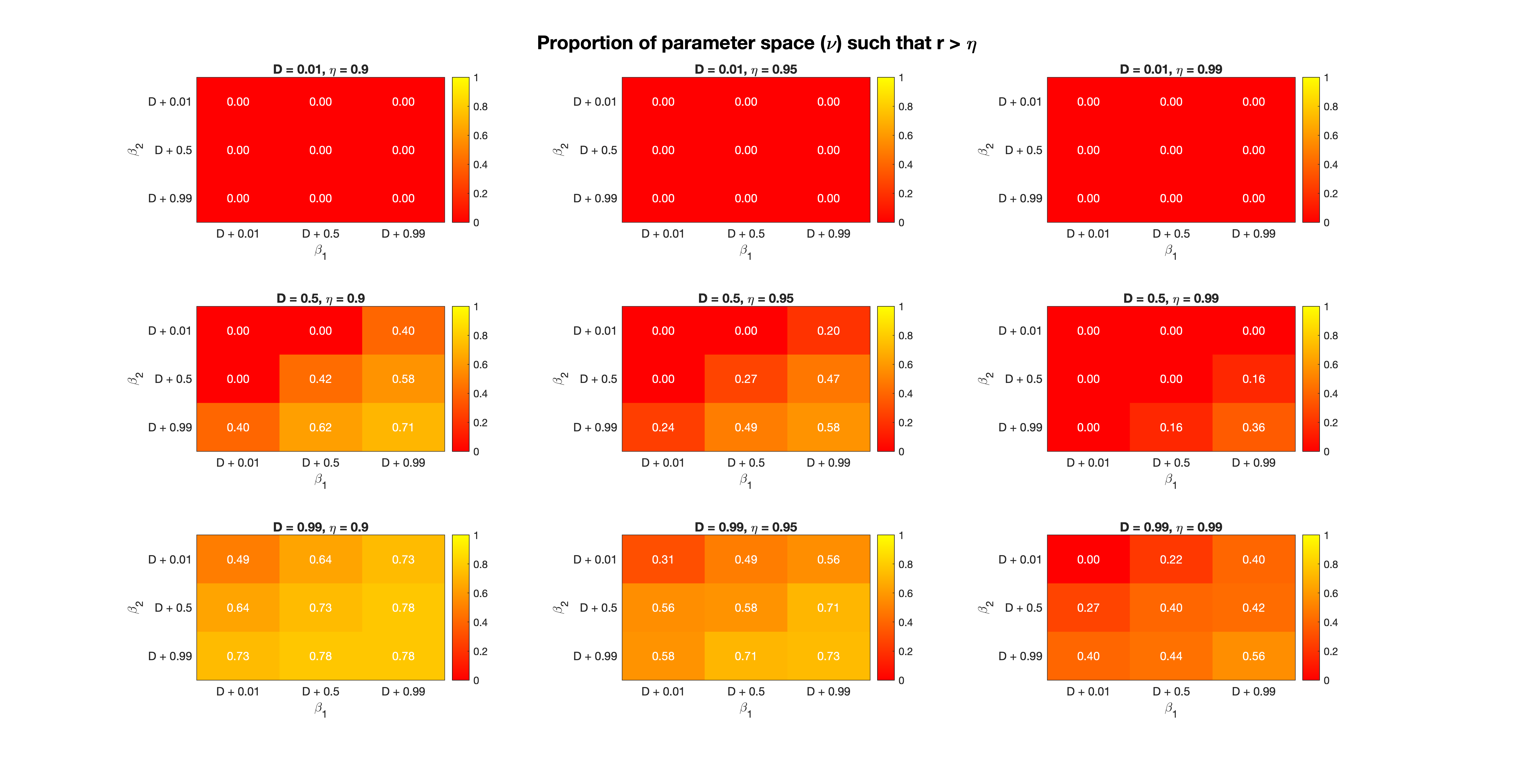}
\caption{\textbf{Exploration of system resilience as a function of model parameters.} A matrix of heatmaps summarizing the simulation study. Each heatmap indicates the value of $\nu$ corresponding to a combination $(\beta_1, \beta_2, D, \eta)$. The range of $\beta_1$ values are on the x-axis and the range of values for $\beta_2$ are on the y-axis for each heatmap. In the matrix comprising the figure, $D$ increases with row number and $\eta$ increases with column number. Since $0 < \nu < 1$, the colorbar for every heatmap ranges from $0$ to $1$.}
\label{three}
\end{figure}

Referring to Fig. \ref{three}, several observations can be made. For instance, low dispersal ($D = 0.01$) guarantees that $\nu = 0$ for $\eta \ge 0.9.$ This means that the desired probability of system resiliency is never $90\%$ or higher if the connectivity between the two patches is low. This makes sense, because the aid of rescue effects in such a system is muted, due to a poor degree of communication between the two populations. Also, $\nu$ decreases as $\eta$ is increased. This is plausible since $\eta$ behaves as a constraint on the parameter space and increasing $\eta$ implies greater specificity. In all cases, if both component Allee effects are weakly strong (i.e. $(\beta_1, \beta_2) = (D + 0.99, D + 0.99)$), then $\nu$ is maximized. This is especially noticeable when $D = 0.99$, regardless of the value of the threshold parameter $\eta$. Lastly, if each heatmap is treated as a matrix, we see that the column (row) of a given matrix is increasing in $\nu$, as the corresponding column (row) indices are increased. This is also reflected within the coloring scheme of each matrix. As expected, $\nu$ is maximal with a value of $0.78$, for $D = 0.99, \eta = 0.9, \beta_1 = D + 0.99,$ and $\beta_2 = D + 0.99.$ \\

\begin{figure}
\centering
\includegraphics[width=\linewidth]{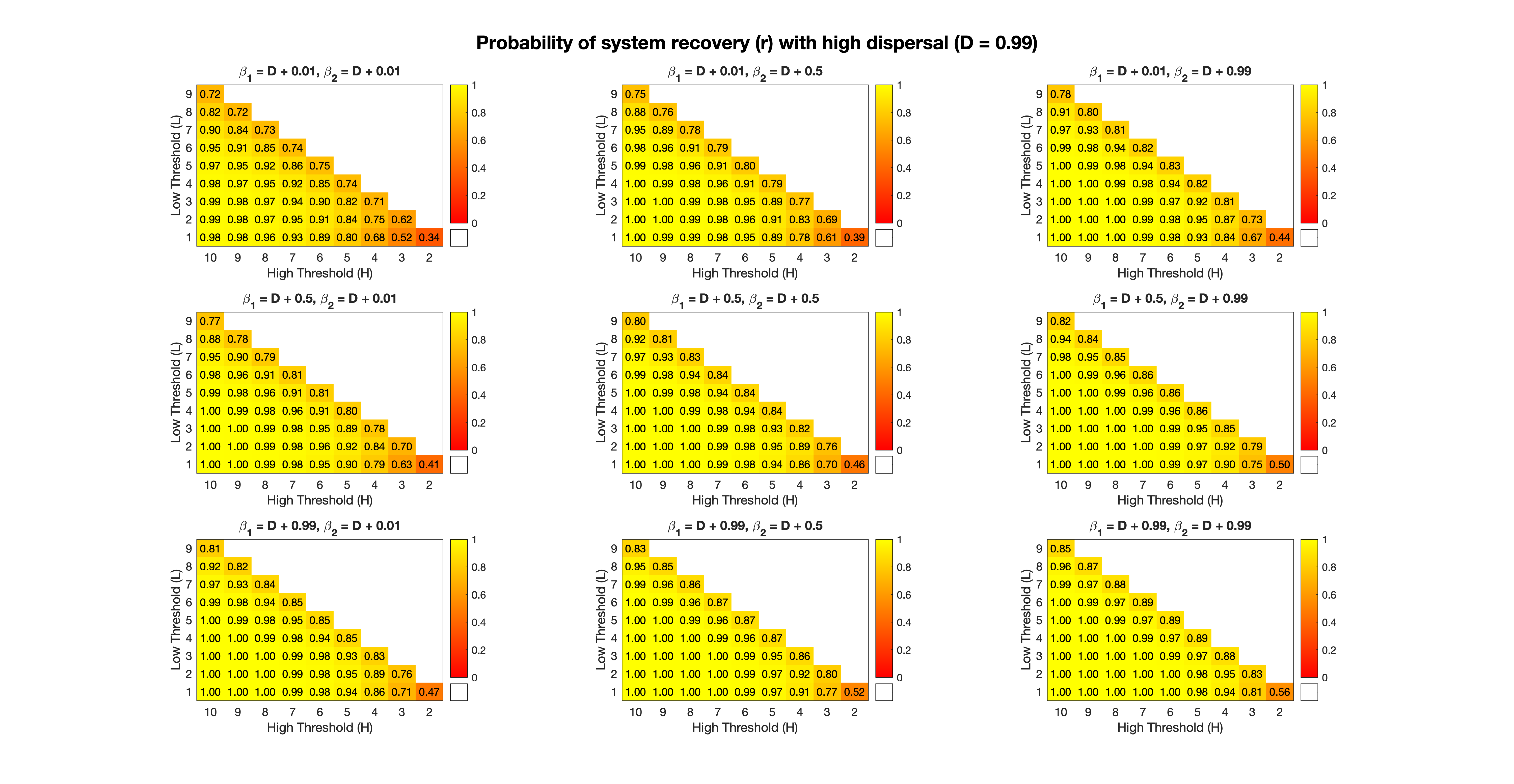}
\caption{\textbf{Likelihood of tipping in the presence of high dispersal.} A representation of the important case of high dispersal corresponding to Fig. \ref{three}. The combinations of $H$ and $L$ are shown explicitly. Note that $H > L$ always. Each heatmap shows the value of $r$ corresponding to a combination $(\beta_1, \beta_2, H, L)$ for $D = 0.99$. In the matrix comprising the figure, $\beta_1$ increases with row number and $\beta_2$ increases with column number. Since $r$ is a probability, each colorbar ranges from $0$ to $1$.}
\label{four}
\end{figure}

In Fig. \ref{four}, we have made explicit the correspondence between Fig. \ref{three} and the combination of high/low thresholds. In concert with the analysis in Fig. \ref{three}, we focused here on the case of high dispersal (i.e. $D = 0.99$). We see that the combination $(H,L) = (2,1)$ is systematically undesirable for the system, as it always results in the lowest value of $r$ for a given pair $(\beta_1,\beta_2)$. However, as the high threshold is increased while the low threshold is kept at the same level, we see that $r$ increases consistently. This is true for any value of the low threshold (i.e. with only $H$ increased). In general, a moderate level of separation between the high and low thresholds is needed to guarantee a desirable outcome for the system. In other words, in order to tip positive change with full confidence ($r \approx 1$), the difference $H - L$ should be sufficiently large. This is sensible in the context of rescue effects, since the high and low thresholds are effectively separated, resulting in a resilience-averaging effect for the system. If the two thresholds were commensurate, the benefit to either population would be reduced due to a lower margin of error between a potential catastrophe or an eventual success. This can be seen clearly in the top-most diagonal of the matrix for each case, as the values of $r$ are low compared to their neighbors, as expected. \\

\begin{figure}
\centering
\includegraphics[width=\linewidth]{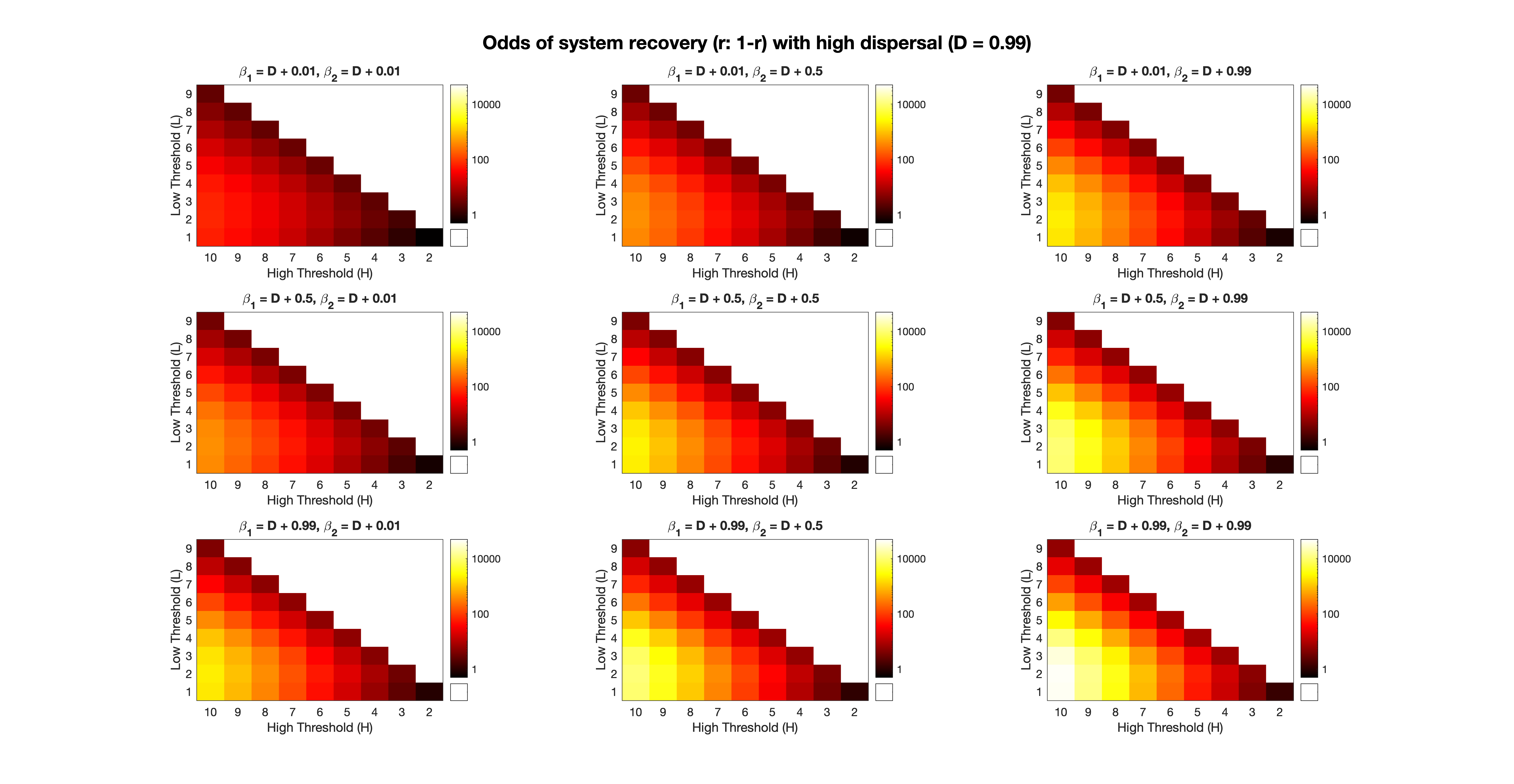}
\caption{\textbf{Odds of tipping with high network connectivity.} This matrix of heatmaps is an analogue to Fig. \ref{four}. The combinations of $H$ and $L$ are shown explicitly in color. Note that $H$ is always larger than $L$. Each heatmap shows the value of $r/(1-r)$, or the odds of system recovery, corresponding to a combination $(\beta_1, \beta_2, H, L)$ for $D = 0.99$. In the matrix comprising the figure, $\beta_1$ increases with row number and $\beta_2$ increases with column number. Each colorbar is on a logarithmic scale.}
\label{five}
\end{figure}

A key issue which relates to the idea of a rescue effect is the relative likelihood of the transition from $HH$ to $LH$, as compared with the transitions $LH$ to $LL$ and $LH$ to $HH$. Referring to Eqs. \eqref{41} - \eqref{43}, we can reason as follows. The likelihood of system resiliency in the form of a rescue effect is given by $P(LH \to HH) = r$, whereas the probability of a total failure or catastrophic collapse is given by $P(LH \to LL) = 1 - r$. Note that \eqref{41} indicates that a partial failure of the system is inevitable (i.e. with probability $1 = r + (1 - r)$). This result follows from the nature of the one-step transition probabilities in the emulator framework. Thus, we can reason that the odds of system recovery is equal to $\frac{r}{1-r}$. As illustrated in Fig. \ref{five}, we see that all of the conclusions observed in Fig. \ref{four} are consistent here as well. For instance, the system is extremely likely to tip favorably in the presence of high dispersal, if both Allee effects are weakly strong. In symbols, $D = 0.99, H = 10, L = 1$, and $\beta_1 = \beta_2 = D + 0.99$. The odds of tipping in this case are nearly $50,000.$ \\

Given this series of observations, there is an interesting and intuitive explanation that generalizes the given setting. Noting that $0 < r < 1$ since $r$ is a probability, we have that $$\dfrac{r}{1-r} = r + r^2 + r^3 + \ldots$$ is a convergent geometric series. Thus, we can write the left-hand side of the equation above as the odds of one network tipping favorably, where the network consists of a system with two patches. The right-hand side, however, can be interpreted as the cumulative probability of all such networks tipping favorably, viz. $$\dfrac{r}{1-r} = \lim\limits_{N \to \infty}\sum\limits_{n = 1}^{N} r^n,$$ where $n$ is the number of two-patch systems in a network. With this analysis, we thus see an instance of a tipping cascade of networks in the form of a domino effect. Note that the propagation of domino dynamics through a network requires strong connectivity, or high dispersal, in our model \citep{lenton-2020}.

\section{Discussion}
Our analysis of a minimal, stochastic system that could have cascading tipping points yields some interesting ecological and mathematical insights. A novel result of this work is that noise-induced tipping can distinguish a catastrophic collapse and a successful recovery from the brink of extinction. This idea is expressed succinctly as the distinction between the tipping of positive change and the failure to adequately address an impending critical transition \citep{lenton-2020}. In our two-population model, we showed that if the system is in the low/high state, a stochastic perturbation due to demographic noise results in either a collapse of both populations (low/low state) or a full recovery to the high/high stable state. This may be a key feature of spatially connected populations with strong Allee effects, and analyzing these dynamics should inform management strategies for similar ecological systems. We find that the population with higher resilience, in the form of a larger Allee threshold, has the capacity to save its counterpart through a rescue effect. This occurs in particular due to the presence of dispersal between the populations. In the presence of noise, patch homogeneity, and strong network connectivity, the system is more likely to exhibit a tipping point cascade \citep{lenton-2020}. \\

Some discussion of systems with noise-induced tipping and their features is needed. In situations where tipping one system increases the probability of another system tipping - for example, melting of the Greenland ice sheet increases the likelihood of failure of the Atlantic Meridional Overturning Circulation (AMOC) - the first system should act as a proxy for the entire network, in terms of a call to action. This is an important point, because by definition, our system can tip without warning, thus precluding detection through generic early-warning signals. However, we have been able to characterize the likelihood of its propensity to tip, either favorably or catastrophically \citep{lenton-2020}. The results from our analysis support the hypothesis that both stochastic and transient dynamics might play a key role in regime shifts and their associated tipping \citep{rocha-2018, hastings-2018}. In particular, domino effects have relatively slow temporal dynamics and larger spatial scales; in our case, the dispersal parameter dominates the time scale in magnitude. In addition, the one-step transition probabilities in our model support the notion that structural dependencies manifest as one-way interactions for the domino effect. This is a salient feature of regime shift couplings. In our system, we note that a partial collapse is inevitable. Most examples of regime shifts exhibiting domino effects concern the Earth system, including monsoon weakening and thermohaline circulation collapse, as well as nutrient transport mechanisms \citep{scheffer, rocha-2018}. \\

Our analysis is based on an individual-based model that describes the dynamics of the system under consideration and this approach should have general applicability. The mean-field description of this model necessarily includes a quadratic term. In other work \citep{abraham-1991, klose-2020}, the proposed tipping element omits the quadratic term, in an attempt to describe a ``dangerous" bifurcation in the form of a cusp catastrophe. This would correspond to zero mortality due to natural causes in our model (i.e. $\gamma = 0$). Similarly, the absence of the quadratic term would indicate that both the carrying capacity and Allee threshold for either population are identically zero \citep{johnson-2018}. Thus, we see that important details of the underlying mechanisms of a tipping element can be lost in a phenomenological description. \\

Ways to generalize our work include looking at more than two populations which would allow the study of how specific forms of connectivity could play a role in tipping cascades, or introducing the possibility of more complex dynamics \citep{strogatz}. Allowing for patch heterogeneity in the two populations by distinguishing their stochastic reaction rates could also produce interesting dynamics. We leave the exploration of these aspects to future work. Given the ubiquity of Allee effects \citep{courchamp} and the generality of our model, we believe that we have uncovered important ecological conclusions that are robust and should apply in a variety of settings. Within the field of landscape ecology, rescue effects can be crucial in both metapopulations and species augmentation efforts. Other potential areas of application include communication systems in network theory and regulatory networks in systems biology.

%
\section*{Conflict of interest}
The authors declare that they have no conflicts of interest.

\bibliographystyle{plainnat}

\begin{thebibliography}{31}
\providecommand{\natexlab}[1]{#1}
\providecommand{\url}[1]{\texttt{#1}}
\expandafter\ifx\csname urlstyle\endcsname\relax
  \providecommand{\doi}[1]{doi: #1}\else
  \providecommand{\doi}{doi: \begingroup \urlstyle{rm}\Url}\fi

\bibitem[Abraham(1991)]{abraham-1991}
Ralph Abraham.
\newblock Cuspoidal nets.
\newblock In \emph{Business Cycles}, pages 56--63. Springer, 1991.

\bibitem[Allen(2010)]{allen-2010}
Linda~JS Allen.
\newblock \emph{An introduction to stochastic processes with applications to
  biology}.
\newblock CRC Press, 2010.

\bibitem[Allen and Allen(2003)]{allen-2003}
Linda~JS Allen and Edward~J Allen.
\newblock A comparison of three different stochastic population models with
  regard to persistence time.
\newblock \emph{Theoretical Population Biology}, 64\penalty0 (4):\penalty0
  439--449, 2003.

\bibitem[Boukal and Berec(2002)]{boukal-2002}
David~S Boukal and Lud{\v{e}}k Berec.
\newblock Single-species models of the allee effect: extinction boundaries, sex
  ratios and mate encounters.
\newblock \emph{Journal of Theoretical Biology}, 218\penalty0 (3):\penalty0
  375--394, 2002.

\bibitem[Chou and D'Orsogna(2014)]{chou-2014}
Tom Chou and Maria~R D'Orsogna.
\newblock First passage problems in biology.
\newblock In \emph{First-passage phenomena and their applications}, pages
  306--345. World Scientific, 2014.

\bibitem[Courchamp et~al.(1999)Courchamp, Clutton-Brock, and
  Grenfell]{courchamp-1999}
Franck Courchamp, Tim Clutton-Brock, and Bryan Grenfell.
\newblock Inverse density dependence and the allee effect.
\newblock \emph{Trends in ecology \& evolution}, 14\penalty0 (10):\penalty0
  405--410, 1999.

\bibitem[Courchamp et~al.(2008)Courchamp, Berec, and Gascoigne]{courchamp}
Franck Courchamp, Ludek Berec, and Joanna Gascoigne.
\newblock \emph{Allee effects in ecology and conservation}.
\newblock Oxford University Press, 2008.

\bibitem[Dai et~al.(2012)Dai, Vorselen, Korolev, and Gore]{dai-2012}
Lei Dai, Daan Vorselen, Kirill~S Korolev, and Jeff Gore.
\newblock Generic indicators for loss of resilience before a tipping point
  leading to population collapse.
\newblock \emph{Science}, 336\penalty0 (6085):\penalty0 1175--1177, 2012.

\bibitem[Drake and Kramer(2011)]{drake-2011}
JM~Drake and AM~Kramer.
\newblock Allee effects.
\newblock \emph{Nature Education Knowledge}, 2011.

\bibitem[Gardiner(2004)]{gardiner}
Crispin~W Gardiner.
\newblock \emph{Handbook of stochastic methods: for physics, chemistry and the
  natural sciences}.
\newblock Springer, 2004.

\bibitem[Gillespie(1976)]{gillespie-1976}
Daniel~T Gillespie.
\newblock A general method for numerically simulating the stochastic time
  evolution of coupled chemical reactions.
\newblock \emph{Journal of computational physics}, 22\penalty0 (4):\penalty0
  403--434, 1976.

\bibitem[Gillespie(1977)]{gillespie-1977}
Daniel~T Gillespie.
\newblock Exact stochastic simulation of coupled chemical reactions.
\newblock \emph{The journal of physical chemistry}, 81\penalty0 (25):\penalty0
  2340--2361, 1977.

\bibitem[Hastings and Gross(2012)]{etc}
Alan Hastings and Louis Gross.
\newblock \emph{Encyclopedia of Theoretical Ecology}.
\newblock Univ of California Press, 2012.

\bibitem[Hastings et~al.(2018)Hastings, Abbott, Cuddington, Francis, Gellner,
  Lai, Morozov, Petrovskii, Scranton, and Zeeman]{hastings-2018}
Alan Hastings, Karen~C Abbott, Kim Cuddington, Tessa Francis, Gabriel Gellner,
  Ying-Cheng Lai, Andrew Morozov, Sergei Petrovskii, Katherine Scranton, and
  Mary~Lou Zeeman.
\newblock Transient phenomena in ecology.
\newblock \emph{Science}, 361\penalty0 (6406), 2018.

\bibitem[Holling(1973)]{holling-1973}
Crawford~S Holling.
\newblock Resilience and stability of ecological systems.
\newblock \emph{Annual review of ecology and systematics}, 4\penalty0
  (1):\penalty0 1--23, 1973.

\bibitem[Johnson and Hastings(2018)]{johnson-2018}
Carter~L Johnson and Alan Hastings.
\newblock Resilience in a two-population system: interactions between allee
  effects and connectivity.
\newblock \emph{Theoretical Ecology}, 11\penalty0 (3):\penalty0 281--289, 2018.

\bibitem[Klose et~al.(2020)Klose, Karle, Winkelmann, and Donges]{klose-2020}
Ann~Kristin Klose, Volker Karle, Ricarda Winkelmann, and Jonathan~F Donges.
\newblock Emergence of cascading dynamics in interacting tipping elements of
  ecology and climate.
\newblock \emph{Royal Society Open Science}, 7\penalty0 (6):\penalty0 200599,
  2020.

\bibitem[Lenton(2020)]{lenton-2020}
Timothy~M Lenton.
\newblock Tipping positive change.
\newblock \emph{Philosophical Transactions of the Royal Society B},
  375\penalty0 (1794):\penalty0 20190123, 2020.

\bibitem[MATLAB(2020)]{MATLAB:2020}
MATLAB.
\newblock \emph{R2020a}.
\newblock The MathWorks Inc., 2020.

\bibitem[M{\'e}ndez et~al.(2019)M{\'e}ndez, Assaf, Mas{\'o}-Puigdellosas,
  Campos, and Horsthemke]{mendez-2019}
Vicen{\c{c}} M{\'e}ndez, Michael Assaf, Axel Mas{\'o}-Puigdellosas, Daniel
  Campos, and Werner Horsthemke.
\newblock Demographic stochasticity and extinction in populations with allee
  effect.
\newblock \emph{Physical Review E}, 99\penalty0 (2):\penalty0 022101, 2019.

\bibitem[Odum and Allee(1954)]{odum-1954}
Howard~T Odum and WC~Allee.
\newblock A note on the stable point of populations showing both intraspecific
  cooperation and disoperation.
\newblock \emph{Ecology}, 35\penalty0 (1):\penalty0 95--97, 1954.

\bibitem[O'Regan(2018)]{oregan-2018}
Suzanne~M O'Regan.
\newblock How noise and coupling influence leading indicators of population
  extinction in a spatially extended ecological system.
\newblock \emph{Journal of biological dynamics}, 12\penalty0 (1):\penalty0
  211--241, 2018.

\bibitem[Polizzi et~al.(2016)Polizzi, Therien, and Beratan]{polizzi-2016}
Nicholas~F Polizzi, Michael~J Therien, and David~N Beratan.
\newblock Mean first-passage times in biology.
\newblock \emph{Israel journal of chemistry}, 56\penalty0 (9-10):\penalty0
  816--824, 2016.

\bibitem[Rocha et~al.(2018)Rocha, Peterson, Bodin, and Levin]{rocha-2018}
Juan~C Rocha, Garry Peterson, {\"O}rjan Bodin, and Simon Levin.
\newblock Cascading regime shifts within and across scales.
\newblock \emph{Science}, 362\penalty0 (6421):\penalty0 1379--1383, 2018.

\bibitem[Scheffer(2009)]{scheffer}
Marten Scheffer.
\newblock \emph{Critical transitions in nature and society}, volume~16.
\newblock Princeton University Press, 2009.

\bibitem[Scheffer et~al.(2012)Scheffer, Carpenter, Lenton, Bascompte, Brock,
  Dakos, Van~de Koppel, Van~de Leemput, Levin, Van~Nes, et~al.]{scheffer-2012}
Marten Scheffer, Stephen~R Carpenter, Timothy~M Lenton, Jordi Bascompte,
  William Brock, Vasilis Dakos, Johan Van~de Koppel, Ingrid~A Van~de Leemput,
  Simon~A Levin, Egbert~H Van~Nes, et~al.
\newblock Anticipating critical transitions.
\newblock \emph{Science}, 338\penalty0 (6105):\penalty0 344--348, 2012.

\bibitem[Strogatz(2001)]{strogatz}
Steven Strogatz.
\newblock \emph{Nonlinear dynamics and chaos: with applications to physics,
  biology, chemistry, and engineering (studies in nonlinearity)}.
\newblock Westview Press, 2001.

\bibitem[van Doorn and Pollett(2013)]{doorn-2013}
Erik~A van Doorn and Philip~K Pollett.
\newblock Quasi-stationary distributions for discrete-state models.
\newblock \emph{European journal of operational research}, 230\penalty0
  (1):\penalty0 1--14, 2013.

\bibitem[Van~Kampen(1992)]{vankampen}
Nicolaas~Godfried Van~Kampen.
\newblock \emph{Stochastic processes in physics and chemistry}.
\newblock Elsevier, 1992.

\bibitem[Volterra(1938)]{volterra-1938}
Vito Volterra.
\newblock Population growth, equilibria, and extinction under specified
  breeding conditions: a development and extension of the theory of the
  logistic curve.
\newblock \emph{Human Biology}, 10\penalty0 (1):\penalty0 1--11, 1938.

\bibitem[Vortkamp et~al.(2020)Vortkamp, Schreiber, Hastings, and
  Hilker]{vortkamp-2020}
Irina Vortkamp, Sebastian~J Schreiber, Alan Hastings, and Frank~M Hilker.
\newblock Multiple attractors and long transients in spatially structured
  populations with an allee effect.
\newblock \emph{Bull Math Biol}, 82\penalty0 (82), 2020.

\end{thebibliography}

\end{document}